\documentclass[12pt]{article}
\usepackage[reqno]{amsmath}
\usepackage{amssymb, theorem, enumerate}
\usepackage{graphicx}
\usepackage{tikz}
\usetikzlibrary{decorations.markings}
\usepackage{array}
\usepackage{float}
\usepackage{wrapfig}
\usepackage{fullpage}

\expandafter\ifx\csname pdfpagebox\endcsname\relax\else
\fi

\theoremstyle{change}
{\theorembodyfont{\slshape}
\newtheorem{theorem}{Theorem.}[section]

\newtheorem{corollary}[theorem]{Corollary.}}

\theorembodyfont{\rmfamily}

\newcommand\cref[1]{Corollary~\ref{cor:#1}}

\def\proof{\noindent{{\sl Proof. }}}
\def\sqr#1#2{{\vbox{\hrule height.#2pt
    \hbox{\vrule width.#2pt height#1pt \kern#1pt
        \vrule width.#2pt}\hrule height.#2pt}}}
\def\eqed{\sqr53}
\def\qed{%
    \ifmmode\eqno\eqed
    \else\nobreak\ \hfill\eqed\medbreak\fi}

\newcommand\cC{{\mathcal C}}

\newcommand\Zy{{\mathbf y}}
\newcommand\Zv{{\mathbf v}}

\newcommand\Zw{{\mathbf w}}

\newcommand\cx{{\mathbb C}}

\newcommand\re{{\mathbb R}}

\newcommand\pmat[1]{\begin{pmatrix} #1 \end{pmatrix}}

\DeclareMathOperator{\col}{col}

\newcommand\ins{D_{h}}
\newcommand\outs{D_{t}}

\title{Cycle Spaces of Digraphs}
\author{Chris Godsil\thanks{Department of Combinatorics \& Optimization, University of Waterloo, Waterloo, ON N2L 3G1.  \texttt{\{cgodsil, kguo\}@uwaterloo.ca} }%
  \and Krystal Guo\footnotemark[1]}

\begin{document}
\maketitle

\begin{abstract}
The cycle space of a graph corresponds to the kernel of an incidence matrix. We investigate an analogous subspace for digraphs. In the case of digraphs of graphs, where every edge is replaced by two oppositely directed arcs, we give a combinatorial description of a basis of such a space. We are motivated by a connection to the transition matrices of discrete-time quantum walks. 
\end{abstract}

\section{Introduction}\label{sec:intro}

The spectrum of a graph is a graph invariant; isomorphic graphs have the same eigenvalues.  The natural follow-up problem is to find classes of graphs in which the spectrum determines the graph, see for example \cite{vDaHa03}. The authors of  \cite{BH, vDaHa03} suggest that the portion of graphs on $n$ vertices which are determined by their spectrum goes to $1$ as $n$ tends to infinity. It is well-known that any two co-parametric strongly regular graphs are cospectral, so it would seem that spectra may not be the right tool to solving Graph Isomorphism in the class of strongly regular graphs. However, there are many proposed algorithms for Graph Isomorphism for the class of strongly regular graphs, which are based on the spectrum of a matrix associated with the graph, see for example \cite{AuGoRoRu07, GFZJ10}. These algorithms are surprisingly successful on small strongly regular graphs. One such proposed routine is that of  Emms et al \cite{EHSW06,ESWH}, which is based on the transition matrix of a discrete-time quantum walk in the arc-reversal model, introduced by Kendon in \cite{Ke06}.  Many classes of strongly regular graphs on were distinguished by this graph invariant and only one pair of counterexamples is known \cite{GGM15}. 

In the course of trying to find an infinite family of counterexamples for the above procedure, the authors were led to the study of a directed version of the cycle space of a graph. We find a decomposition of the vector space indexed by the arcs of a graph, with the goal of using it to diagonalize matrices indexed by the arcs of graph.  The relation between discrete-time quantum walks and Ihara zeta function was observed in \cite{raewh11}. The positive support of the transition matrix is related to the Bass-Hashimoto edge adjacency operator. 

We proceed with a few basic definitions. The \textsl{digraph of a graph $G$} is obtained from $G$ by replacing each edge $\{x,y\}$ with directed edges $xy$ and $yx$.  In this article, we study a related kernel $\ker(\outs) \cap \ker(\ins)$, for $X$ a digraph of a graph $G$. We are interested in its dimension and descriptions of the combinatorial objects in this space.

Let $X$ be a digraph without loops or parallel arcs. Let $A$ the adjacency matrix of $X$. We consider the following incidence matrices of $X$, both with rows indexed by the vertices of $X$ and columns indexed by the arcs of $X$:
\[ (\ins)_{i,j} = \begin{cases} 1, &\text{if } i \text{ is the head of arc }j; \\ 
0, &\text{otherwise,}\end{cases}
\]
and 
\[ (\outs)_{i,j}  = \begin{cases} 1, &\text{if } i \text{ is the tail of arc }j ;\\ 
0, &\text{otherwise.}\end{cases}
\]
Let $B(X) = \outs + \ins$ and $N(X) = \outs - \ins$. We will write $N$ and $B$ for simplicity, when the choice of digraph is clear. 

In this paper, we find the following main result. 

\begin{theorem} \label{thm:rk} For $X$ the digraph of a graph $G$ with $m$ edges, the subspace $\ker(\outs) \cap \ker(\ins)$ of $\re^{2m}$ has dimension $2m - 2n +b +c$, where $b$ is the number of bipartite components of $X$ and $c$ is the number of components of $X$. \end{theorem}

In the case that $G$ is a bipartite graph, we can give an explicit basis of $\ker(\outs) \cap \ker(\ins)$ in terms of the cycle space of $G$.

\section{Preliminaries}

Given a graph $X$, there are two incidence matrices of $X$, which are commonly studied. We will follow the notation and definitions in \cite{GR}. For an orientation $X^{\sigma}$ of $X$, they are given as follows:
\[
B(X) = D_h(X^{\sigma}) + D_t(X^{\sigma})
\]
and
\[
N(X^{\sigma}) = D_t(X^{\sigma}) - D_h(X^{\sigma}),
\]
where $D_t(X^{\sigma})$ and $D_h(X^{\sigma})$ are the tails and heads incidence matrices of $X^{\sigma}$. 

\begin{theorem}[Theorem 8.2.1, \cite{GR}] For a graph $X$ on $n$ vertices with $b$ bipartite components, the incidence matrix $B(X)$ has rank equal to $n-b$. \end{theorem}
 
It is apparent from the definition that any choice of orientation gives the same incidence matrix $B(X)$. The following theorem implies the choice of orientation does not affect that rank of $N(X^{\sigma})$. 
 
 \begin{theorem}[Theorem 8.3.1, \cite{GR}] For a graph $X$ on $n$ vertices with $c$  components and an orientation $X^{\sigma}$ of $X$, the incidence matrix $N(X^{\sigma})$ has rank equal to $n-c$. \end{theorem}
 
The kernel of $N(X^{\sigma})$ is called the \textsl{flow space} or \textsl{cycle space} and has a combinatorial description in terms of the cycles of the graph. Let $C$ be a cycle in $X$; that is $C$ is a set of directed edges of $X$ such that if $uv \in C$ then $vu \notin C$ and $C$ induces a cycle in $X$, when the directions are forgotten. With respect to an orientation $X^{\sigma}$ of $X$, the \textsl{signed characteristic vector} of $C$ in $\cx^{m}$ is as follows:
\[
(\Zv_C)_{uv} = \begin{cases}1, &\text{ if } uv \in C;\\
-1, &\text{ if } vu \in C; \\
0, &\text{otherwise}. \end{cases} 
\]
Note that for each directed edge $uv$ of $C$, exactly one of $uv$ and $vu$ appear as a directed edge of $X^{\sigma}$. 

 \begin{theorem}[Corollary 14.2.3, \cite{GR}] For a graph $X$ with an orientation $X^{\sigma}$ of $X$, the flow space $\ker(N(X^{\sigma}))$ is spanned by the signed characteristic vectors of its cycles. \end{theorem}

\section{Main result}

\begin{theorem} \label{thm:rk} For $X$ the digraph of a graph $G$ on $n$ vertices and $m$ edges, the subspace $\ker(\outs(X)) \cap \ker(\ins(X))$ of $\re^{2m}$ is
\[
\left\{ H\pmat{\Zv \\ \Zw} \mid \Zv \in \ker(B(X^{\sigma})),\, \Zw \in \ker(N(X^{\sigma})) \right\}
\] 
where 
\[
H = \frac{1}{\sqrt{2}} \pmat{I & I \\I & -I},
\]
and thus has dimension $2m -2n +b+c$, where $b$ is the number of bipartite components of $X$ and $c$ is the number of components of $X$. \end{theorem}

\proof First observe the following:
\[ \ker(\outs(X)) \cap \ker(\ins(X))= \ker\pmat{\outs(X) \\ \ins(X)}.
\]

We consider $X^{\sigma}$, any orientation of $G$. Let the edges of $X^{\sigma}$ be $\{e_1, \ldots, e_m\}$. We order the edges of $X$ as $\{e_1, \ldots, e_m, \bar{e}_1, \ldots, \bar{e}_m\}$ where if $e_i = uv$ then $\bar{e}_e = vu$. We have that 
\[
\outs(X) = \pmat{\outs(X^{\sigma}) & \ins(X^{\sigma})} \text{ and } \ins(X) =\pmat{\ins(X^{\sigma}) & \outs(X^{\sigma})}.
\]
Thus
\[
\pmat{\outs(X) \\ \ins(X)} = \pmat{\outs(X^{\sigma}) & \ins(X^{\sigma})\\ \ins(X^{\sigma}) & \outs(X^{\sigma})}.
\]
We see that
\[
H \pmat{\outs(X^{\sigma}) & \ins(X^{\sigma})\\ \ins(X^{\sigma}) & \outs(X^{\sigma})} H = \pmat{\outs(X^{\sigma}) + \ins(X^{\sigma}) & 0 \\ 0& \ins(X^{\sigma}) - \outs(X^{\sigma})} = \pmat{B(X^{\sigma}) & 0 \\ 0& N(X^{\sigma})}
\]
where 
\[
H = \frac{1}{\sqrt{2}} \pmat{I & I \\I & -I}
\]
and $H^T = H^{-1} = H$.  We have that
\[
\ker\pmat{B(X^{\sigma}) & 0 \\ 0& N(X^{\sigma})} = \left\{ \pmat{\Zv \\ \Zw} \mid \Zv \in \ker(B(X^{\sigma})),\, \Zw \in \ker(N(X^{\sigma})) \right\}.
\] 
For $\Zv \in \ker(B(X^{\sigma}))$ and  $\Zw \in \ker(N(X^{\sigma}))$, we consider the vector $H \pmat{\Zv \\ \Zw}$:
\[
\begin{split}
\pmat{\outs(X) \\ \ins(X)} H \pmat{\Zv \\ \Zw} &= H \pmat{B(X^{\sigma}) & 0 \\ 0& N(X^{\sigma})} H H \pmat{\Zv \\ \Zw} \\
&= H \pmat{B(X^{\sigma}) & 0 \\ 0& N(X^{\sigma})} \pmat{\Zv \\ \Zw} \\
&= 0
\end{split}
\]
Thus,
\[
\ker\pmat{\outs(X) \\ \ins(X)} = \left\{ H\pmat{\Zv \\ \Zw} \mid \Zv \in \ker(B(X^{\sigma})),\, \Zw \in \ker(N(X^{\sigma})) \right\}
\]
and
\[
\dim \ker(\outs(X)) \cap \ker(\ins(X)) = 2m - 2n + c + b.
\]
as claimed.\qed 

When $G$ is bipartite, there is a natural choice of an orientation, which, together with the matrix of similarity, gives an explicit basis for $L$. 

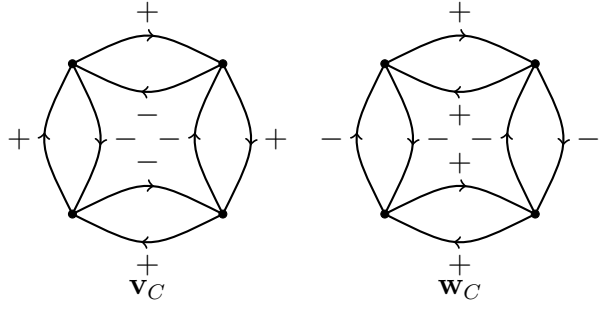
\begin{figure}[ht]
\begin{center}
\begin{tikzpicture}[decoration={
markings,
mark=at position .54 with {\arrow[black,thick]{>};}}
]
\filldraw[black]  (-1,-1) circle (1.5pt)
                           (1,-1) circle (1.5pt)
                           (-1,1) circle (1.5pt)
                           (1,1) circle (1.5pt);
\draw[postaction={decorate},thick] (-1,1) .. controls (0, 1.5)  .. (1,1);
\draw[postaction={decorate},thick] (1,1) .. controls (0, 0.5)  .. (-1,1);
\draw[postaction={decorate},thick] (-1,-1) .. controls (0, -0.5)  .. (1,-1);
\draw[postaction={decorate},thick] (1,-1) .. controls (0, -1.5)  .. (-1,-1);
\draw[postaction={decorate},thick] (1,1) .. controls (1.5, 0)  .. (1,-1);
\draw[postaction={decorate},thick] (1,-1) .. controls (0.5, 0)  .. (1,1);
\draw[postaction={decorate},thick] (-1,1) .. controls (-0.5, 0)  .. (-1,-1);
\draw[postaction={decorate},thick] (-1,-1) .. controls (-1.5, 0)  .. (-1,1);
\draw (0,1.4) node[anchor=south] {$+$}
           (0,0.6) node[anchor= north] {$-$}
           (0,-0.6) node[anchor=south] {$-$}
           (0,-1.4) node[anchor=north] {$+$}
           (-1.4,0) node[anchor=east] {$+$}
           (-0.6,0) node[anchor=west] {$-$}
           (0.6,0) node[anchor=east] {$-$}
           (1.4,0) node[anchor=west] {$+$}
           (0, -2) node[anchor=center] {$\Zv_C$};
\end{tikzpicture} 
\begin{tikzpicture}[decoration={
markings,
mark=at position .54 with {\arrow[black,thick]{>};}}
]
\filldraw[black]  (-1,-1) circle (1.5pt)
                           (1,-1) circle (1.5pt)
                           (-1,1) circle (1.5pt)
                           (1,1) circle (1.5pt);
\draw[postaction={decorate},thick] (-1,1) .. controls (0, 1.5)  .. (1,1);
\draw[postaction={decorate},thick] (1,1) .. controls (0, 0.5)  .. (-1,1);
\draw[postaction={decorate},thick] (-1,-1) .. controls (0, -0.5)  .. (1,-1);
\draw[postaction={decorate},thick] (1,-1) .. controls (0, -1.5)  .. (-1,-1);
\draw[postaction={decorate},thick] (1,1) .. controls (1.5, 0)  .. (1,-1);
\draw[postaction={decorate},thick] (1,-1) .. controls (0.5, 0)  .. (1,1);
\draw[postaction={decorate},thick] (-1,1) .. controls (-0.5, 0)  .. (-1,-1);
\draw[postaction={decorate},thick] (-1,-1) .. controls (-1.5, 0)  .. (-1,1);
\draw (0,1.4) node[anchor=south] {$+$}
           (0,0.6) node[anchor= north] {$+$}
           (0,-0.6) node[anchor=south] {$+$}
           (0,-1.4) node[anchor=north] {$+$}
           (-1.4,0) node[anchor=east] {$-$}
           (-0.6,0) node[anchor=west] {$-$}
           (0.6,0) node[anchor=east] {$-$}
           (1.4,0) node[anchor=west] {$-$}
           (0, -2) node[anchor=center] {$\Zw_C$};
\end{tikzpicture} 
\end{center}
\caption{Vectors  $\Zy_C$ and $\Zw_C$ for a $4$-cycle. \label{fig:yw}}
\end{figure}

With respect to $X^{\sigma}$, we define two characteristic vectors $\Zy_C$ and $\Zw_C$ of $C$ in $\cx^{2m}$ as follows:
\[
(\Zy_C)_{uv} = \begin{cases} \Zv_{uv}, &\text{ if } uv \in E(X^{\sigma});\\
-\Zv_{uv}, &\text{otherwise} \end{cases} 
\] and
\[
(\Zw_C)_{uv} = \begin{cases} \Zv_{uv}, &\text{ if } uv \in E(X^{\sigma});\\
\Zv_{uv}, &\text{otherwise}. \end{cases} 
\]
See Figure \ref{fig:yw}.

\begin{theorem}Let $X$ be a bipartite graph, $X^{\sigma}$ be an orientation of $X$, and $\cC$ be a cycle basis of $X$. The subspace $\dim \ker(\outs(X)) \cap \ker(\ins(X))$ of $\cx^{2m}$ has the following basis:
\[
\{\Zy_C, \Zw_C  \mid C \in \cC \}.
\]  \end{theorem}

\proof We take $X^{\sigma}$ to be the orientation of $X$ with bipartition $(Y,Z)$ to have all directed edges of $X$ with tails in $Y$ and heads in $Z$. Let $R$ be the $n \times n$ diagonal matrix with entries as follows:
\[
R_{v,v} = \begin{cases}1, &\text{ if } v \in Y; \\
 -1, &\text{ if } v \in Z. \end{cases}
\]
We see that $RB(X^{\sigma}) = N(X^{\sigma}) $.  

Let $\Zv \in \ker(N(X^{\sigma}))$. Since $N(X^{\sigma}) \Zv = 0$, we have that $RB(X^{\sigma}) \Zv = 0$. Since $R$ is an invertible matrix, we have that $\Zv \in \ker B(X^{\sigma})$ and so $ \ker(N(X^{\sigma}))\subseteq\ker(B(X^{\sigma}))$. Since $R^2 = I$, we also have $RN(X^{\sigma}) = B(X^{\sigma}) $ and can similarly show that $ \ker(B(X^{\sigma}))\subseteq\ker(N(X^{\sigma}))$, and so $ \ker(N(X^{\sigma}))=\ker( B(X^{\sigma}))$. 

Theorem \ref{thm:rk} gives us that 
\[
\ker\pmat{\outs(X) \\ \ins(X)} = \left\{ H\pmat{\Zv \\ \Zw} \mid \Zv ,\Zw \in \ker(N(X^{\sigma})) \right\} . 
\] 
Observe that 
\[
\left\{ H\pmat{\Zv \\ 0}, H\pmat{0 \\ \Zv} \mid \Zv \in \ker(N(X^{\sigma})) \right\}
\]
is an independent set of vectors  in $\ker\pmat{\outs(X) \\ \ins(X)} $ whose cardinality is equal to the dimension of $\ker\pmat{\outs(X) \\ \ins(X)} $ and is hence a basis. 
Consider $\Zv \in \ker(N(X^{\sigma}))$. We see that 
\[
H \pmat{\Zv \\ 0 } = \frac{1}{\sqrt{2}} \pmat{\Zv \\ \Zv }
\]
and 
\[
H \pmat{0 \\ \Zv } = \frac{1}{\sqrt{2}} \pmat{\Zv \\ -\Zv }.
\]
If $\Zv$ is the characteristic vector of a cycle $C$ in $X^{\sigma}$, then $\pmat{\Zv \\ \Zv } = \Zw_C$ and  $\pmat{\Zv \\ -\Zv } = \Zy_C$.
  \qed 

\section{Relation to quantum walks} 

A \textsl{discrete quantum walk} is a process on a graph $G$ governed by a unitary matrix, $U$, which is called the \textsl{transition matrix.} For $uv$ and $wx$ arcs in the digraph of $X$, the transition matrix is defined to be:
\[
 U_{wx,uv} = \begin{cases} \frac{2}{d(v)} &\text{ if } v=w \text{ and } u \neq x ,\\
\frac{2}{d(v)} -1  & \text{ if } v=w \text{ and } u = x, \\
0 &\text{ otherwise.} \end{cases}
\]

We may write this in terms of the incidence matrices of $X$, the digraph of $G$. To describe the quantum walk, we need one more matrix: let $P$ be a permutation matrix with row and columns indexed by the arcs of $D$ such that,
\[P_{wx,uv}  = \begin{cases} 1 &\text{if } x=u \text{ is the tail of arc } w=v \\ 
0 &\text{otherwise.}\end{cases}
\]
Then, we see that $\ins\outs^{T} = A(G)$, the adjacency matrix of $G$, and
\[
 (\outs^T\ins)_{wx,uv} = \begin{cases} 1 &\text{if } v=w,\\
0 &\text{otherwise.} \end{cases}
\]
If $G$ is regular with valency $k$, we have that \[U = \frac{2}{k} \outs^T\ins -  P.\] It is important to note that $U$ is a $nk \times nk$ unitary matrix. In addition, from the definitions of $\ins$, $\outs$ and $P$, we easily see the following:
\[
\ins \ins^T = k I, \, \outs \outs^T = kI, \, \ins P = \outs, \, \text{ and } \outs P = \ins.
\]

In \cite{GG11}, we prove that one can diagonalize $U$ and related matrices over $\cx^{kn}$ by decomposing $\cx^{kn} = K \oplus L$, where $K = \col(\ins) \oplus \col(\outs)$ and $L = \ker(\outs) \cap \ker(\ins)$.

\section{Further applications}

The Bass-Hashimoto edge adjacency matrix $T(G)$ of a graph $G$ is a matrix indexed by end arcs of $G$ such that 
\[
T(G)_{uv,wx} = \begin{cases} 
1, &\text{ if } v = w \text{ and } u \neq x;\\
0, &\text{ otherwise.}
\end{cases}
\]
The Bass-Hashimoto edge adjacency matrix has been studied in the context of the Ihara zeta function of graphs. [add citations here] Observe that we can write $T(G)$ in terms of incidence matrices as follows:
\[
T(G) = D_h^T D_t - P,
\]
where $P$ is the permutation matrix taking each arc to the reverse arc. When the context is clear, we will write $T$ for $T(G)$. 

An eigenvalue is said to be \textsl{semi-simple} if its algebraic and geometric multiplicities are equal. A matrix is \text{semi-simple} if all of its eigenvalues are semi-simple. In \cite{CoKo15}, the authors find that $T$ is semi-simple over $\ker(\outs(X)) \cap \ker(\ins(X))$ but can fail to be semi-simple in general. In particular, they show that $T$ is not semi-simple if $T$ has a vertex of degree $1$ and they ask if the presence of vertices of degree $1$ are the only obstructions to simplicity. Here, we give an answer in the negative by computation, but find that $T$ is semi-simple for regular graphs of degree at least $2$. 

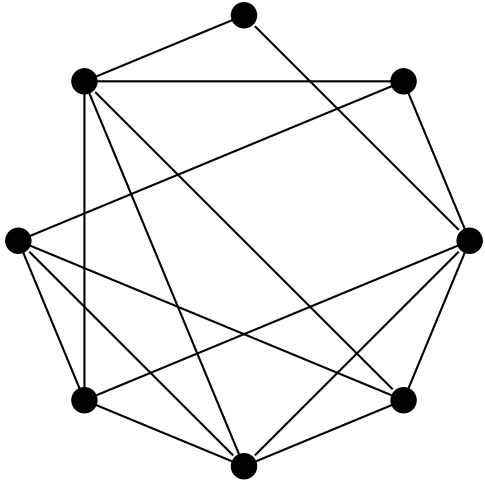
\begin{figure}[ht]
\begin{center}
\begin{tikzpicture}
\foreach \a in {0,1,2,3,4,5,6,7}{
\fill (\a*360/8: 3cm) circle (5pt) node(\a){};
}

\draw[thick] (0) -- (1);
\draw[thick] (0) -- (2);
\draw[thick] (0) -- (5);
\draw[thick] (0) -- (6);
\draw[thick] (0) -- (7);
\draw[thick] (1) -- (3);
\draw[thick] (1) -- (4);
\draw[thick] (2) -- (3);
\draw[thick] (3) -- (5);
\draw[thick] (3) -- (6);
\draw[thick] (3) -- (7);
\draw[thick] (4) -- (5);
\draw[thick] (4) -- (6);
\draw[thick] (4) -- (7);
\draw[thick] (5) -- (6);
\draw[thick] (6) -- (7);
\end{tikzpicture} 
\end{center}
\caption{A graph on $8$ vertices with two non-semi-simple eigenvalues and no vertex of degree $1$. \label{fig:ex8}}
\end{figure}

We find, by a computation using Sage \cite{sage}, that for graphs on $n$ vertices where $n = 1,\ldots, 6$ the only graphs for which the Bass-Hashimoto edge adjacency matrix has a non-semi-simple eigenvalue have a vertex of degree $1$. However, the statement is false for graphs on $7$ vertices; there are $2$ graphs which have $(x^2 + x + 2)^2$ as a factor of the minimal polynomial pf $T$ and so the roots of $(x^2 + x + 2)^2$ are  non-semi-simple eigenvalues. For graphs on $8$ vertices, there are $52$ graphs which have a non-semi-simple eigenvalue and which contain no vertex of degree $1$. Of these graphs $22$ have $x^2+2$ as a repeated root of the minimal polynomial of $T$ and $30$ have $x^2 + x + 2$ as a repeated root of the minimal polynomial of $T$. Figure \ref{fig:ex8} shows an example of one such graph; it is connected with a connected complement, has no vertex of degree $1$, and its  minimal polynomial of $T$ is 
\[
(x - 1)(x + 1) (x^{2} + 3) (x^{2} + 4)  (x^{2} + 2)^{2} (x^{7} + x^{6} - 2x^{5} - 14x^{4} - 39x^{3} - 59x^{2} - 72x - 72).
\]

In relation to discrete-time quantum walks, the authors of \cite{EHSW06} study the positive support of the transition matrix $U$; they study the matrix $S^+(U)$ whose $(i,j)$ entry is $1$ whenever the $(i,j)$ entry of $U$ is positive, and is $0$ otherwise. In terms of incidence matrices, we can write this matrix as follows: 
\[
S^+(U) = D_t^T D_h - P.
\]
Observe that
\[
PT(G)P = PD_h^T D_tP - P^3 = D_t^T D_h - P
\]
which implies that $T(G)$ is similar to $S^+(U)$ via $P$ and hence $T$ and $S^+(U)$ have the same eigenvalues and, further, an eigenvalue $\lambda$ of $T$ is semi-simple if and only if $\lambda$ is a semi-simple eigenvalue of $S^+(U)$. 

Theorem 3.1 of \cite{GG11} finds that for $k$-regular graphs on $n$ vertices with $k \geq 2$, $S^+(U)$ is diagonalizable over $\cx^{kn}$ and find the eigenvalues by showing that $L$ and $K$ are $S^+(U)$-invariant and diagonalizing over each space separately. This gives the immediately corollary.

\begin{corollary} The Bass-Hashimoto edge adjacency matrix $T$ is semi-simple for all $k$-regular graphs with $k\geq 2$. 
\end{corollary}


\begin{thebibliography}{10}

\bibitem{AuGoRoRu07}
Koenraad Audenaert, Chris Godsil, Gordon Royle, and Terry Rudolph.
\newblock Symmetric squares of graphs.
\newblock {\em J. Combin. Theory Ser. B}, 97(1):74--90, 2007.

\bibitem{BH}
Andries~E. Brouwer and Willem~H. Haemers.
\newblock {\em Spectra of graphs}.
\newblock Universitext. Springer, New York, 2012.

\bibitem{CoKo15}
G.~{Cornelissen} and J.~{Kool}.
\newblock {Edge reconstruction of the Ihara zeta function}.
\newblock {\em ArXiv e-prints}, July 2015.

\bibitem{EHSW06}
David Emms, Edwin~R. Hancock, Simone Severini, and Richard~C. Wilson.
\newblock A matrix representation of graphs and its spectrum as a graph
  invariant.
\newblock {\em Electr. J. Comb.}, 13(1), 2006.

\bibitem{ESWH}
David Emms, Simone Severini, Richard~C. Wilson, and Edwin~R. Hancock.
\newblock Coined quantum walks lift the cospectrality of graphs and trees.
\newblock {\em Pattern Recognition}, 42(9):1988--2002, 2009.

\bibitem{GFZJ10}
John~King Gamble, Mark Friesen, Dong Zhou, Robert Joynt, and S.~N. Coppersmith.
\newblock Two-particle quantum walks applied to the graph isomorphism problem.
\newblock {\em Phys. Rev. A}, 81(5):052313, May 2010.

\bibitem{GGM15}
C.~{Godsil}, K.~{Guo}, and T.~G.~J. {Myklebust}.
\newblock {Quantum walks on generalized quadrangles}.
\newblock {\em ArXiv e-prints}, November 2015.

\bibitem{GR}
C.~Godsil and G.~Royle.
\newblock {\em Algebraic Graph Theory}.
\newblock Springer-Verlag, New York, 2001.

\bibitem{GG11}
Chris Godsil and Krystal Guo.
\newblock Quantum walks on regular graphs and eigenvalues.
\newblock {\em Electr. J. Comb.}, 18(1):165, 2011.

\bibitem{Ke06}
Viv Kendon.
\newblock Quantum walks on general graphs.
\newblock {\em International Journal of Quantum Information}, 04(05):791--805,
  2006.

\bibitem{raewh11}
Peng Ren, Tatjana Aleksi{\'c}, David Emms, Richard Wilson, and Edwin Hancock.
\newblock Quantum walks, {I}hara zeta functions and cospectrality in regular
  graphs.
\newblock {\em Quantum Information Processing}, 10:405--417, 2011.
\newblock 10.1007/s11128-010-0205-y.

\bibitem{sage}
W.\thinspace{}A. Stein et~al.
\newblock {\em {S}age {M}athematics {S}oftware ({V}ersion 6.1.1)}.
\newblock The Sage Development Team, 2014.
\newblock {\tt http://www.sagemath.org}.

\bibitem{vDaHa03}
Edwin~R. van Dam and Willem~H. Haemers.
\newblock Which graphs are determined by their spectrum?
\newblock {\em Linear Algebra and its Applications}, 373:241 -- 272, 2003.

\end{thebibliography}

\end{document}